\documentclass[12pt,reqno]{article}

\usepackage[english]{babel}
\usepackage{amsmath,amsfonts,amssymb,amsthm}
\usepackage{mathrsfs}
\usepackage{graphicx}
\usepackage{bbm}
\usepackage{indentfirst}
\usepackage{graphicx}

\theoremstyle{definition}

\numberwithin{equation}{section}

\title{
{\bf\Large Remarks on a rumor propagation model}}
\author{{\bf\large Alberto Ragagnin}
\vspace{1mm}\\
{\it\small University of Udine}\\
{\it\small Department of Mathematics, Computer Science and Physics}\\
{\it\small Via delle Scienze 206, 33100 Udine, Italy}\\
{\it\small e-mail: alberto.ragagnin.85@gmail.com}}

\date{}

\begin{document}

\maketitle

\begin{abstract}
This short note contains a few comments and corrections about
some recent models for the spread of rumors in a population. We
consider a system of ordinary differential equations which describes
the evolution of Ignorant-Spreaders-Stiflers in time. State of the art
of analytical understanding of those equations is based on studying
asymptotic solutions of the rumor spreading equations. In this work
we find a First Integral of these differential equations. We qualitatively
discuss the evolution of the system in the light of those new more
precise solutions.
\end{abstract}
\maketitle

\section{Introduction}\label{sec-1}

The study of the propagation of rumors in a population has become
a research topic of increasing interest in the recent years.
Motivations for these investigations come from different
perspectives, such as social sciences, economy, informatics and
military interests.

The first rumor propagation models considered in the literature
have been adapted from the famous $SIR$ model by Kermack-McKendrik,
whose history is well described in \cite{Ba-11}.
However, in modelling rumor
propagation, the mechanism differs basically from that governing the spread of an epidemic.
Interesting pioneering works were made by
Daley and Kendall \cite{DaKe-65} who proposed a
stochastic, random-walk, model and by
Maki and Thompson \cite{MaTh-73} who treat a deterministic discrete model.
In \cite{MaTh-73} the authors did not write down explicitly the
assumptions for a continuous-time model, but left this task as
an exercise for the reader \cite[Ch. 9, p.388]{MaTh-73}. From this point
of view, the
sources that we have found in literature are somehow a little confusing,
since usually one refers to as Maki-Thompson models also models based on
differential equations.

In this
paper we focus our attention on two differential equation models
investigated by J.R. Piqueira in the recent article \cite{Pi-10}
and S. Belen and C. Pearce in \cite{BePe-04}, respectively.

The set of population could be partitioned into three subsets of
sub-populations: the \textit{$I$-Ignorant}, namely the individuals
who ignore the rumor (who play the same role as the susceptible of the SIR model),
the \textit{$S$-Spreaders} who disseminate the rumor
(and play the same role as the infected of the SIR model),
and the
\textit{$R$-Stiflers} who do not spread the rumor after receiving it
(who play a similar role as the recovered of the SIR model).

\section{Remarks on Piqueira's Model}\label{sec-2}

We have obtained a different result and a different conclusion
for the dynamical system proposed in the original paper \cite{Pi-10} by J.R. Piqueira.
Let us briefly recall Piqueira's model.
The functions $I(t),$ $S(t)$ and $R(t)$ are continuously differentiable and
represent the number of individuals of the three sub-populations
at the time $t$. We suppose that along all the time interval in which we study the model,
the total population is constant, that is
\begin{equation}\label{cond0}
I(t) + S(t) + R(t) = N, \quad\forall\, t.
\end{equation}
It will be not restrictive to suppose $N=1$.
The dynamics of the triplet $(I(t),S(t),R(t))$
is supposed to be governed by the nonlinear ODE:

\begin{equation}\label{eq-1}
\begin{cases}
I'=-\rho_2 \mu IS\\
S'= \rho_2 \mu IS-\rho_1  \mu S(S+R) \\
R'=  \rho_1 \mu S(S+R)\\
\end{cases}
\end{equation}
where the positive parameters $\rho_1,$ $\rho_2$ and $\mu$ are assigned as follows:
$\rho_1$ is the probability that a Spreader meets another Spreader causing their silencing,
$\rho_2$ is the probability that an Ignorant converts into a Spreader after heard the rumor and
$\mu$ is the average number of contacts for every individual
(see \cite{Pi-10}).
According to \cite{Pi-10}, this system is inspired by the Daley-Kendall model, following also
\cite{Mo-04}.
Concerning equation \eqref{eq-1}, we notice that the parameter $\mu$ is
practically useless from the point of view of the \textit{qualitative analysis}
and, therefore, we could omit it  by posing  $\mu=1$. Accordingly, system
\eqref{eq-1} takes the form

\begin{equation}\label{eq-2}
\begin{cases}
I'=-\rho_2 IS\\
S'=  S( \rho_2 I-\rho_1  (S+R)) \\
R'=  \rho_1 S(S+R)\\
\end{cases}
\end{equation}
\subsection{Comparison with Asymptotic Stability}
From now on, we will focus our attention to the qualitative study of the trajectories
for which numerical simulations show a behavior like that described in Figure $1$, below.

\begin{figure}[h]
\centering
\includegraphics[width=0.7\textwidth]{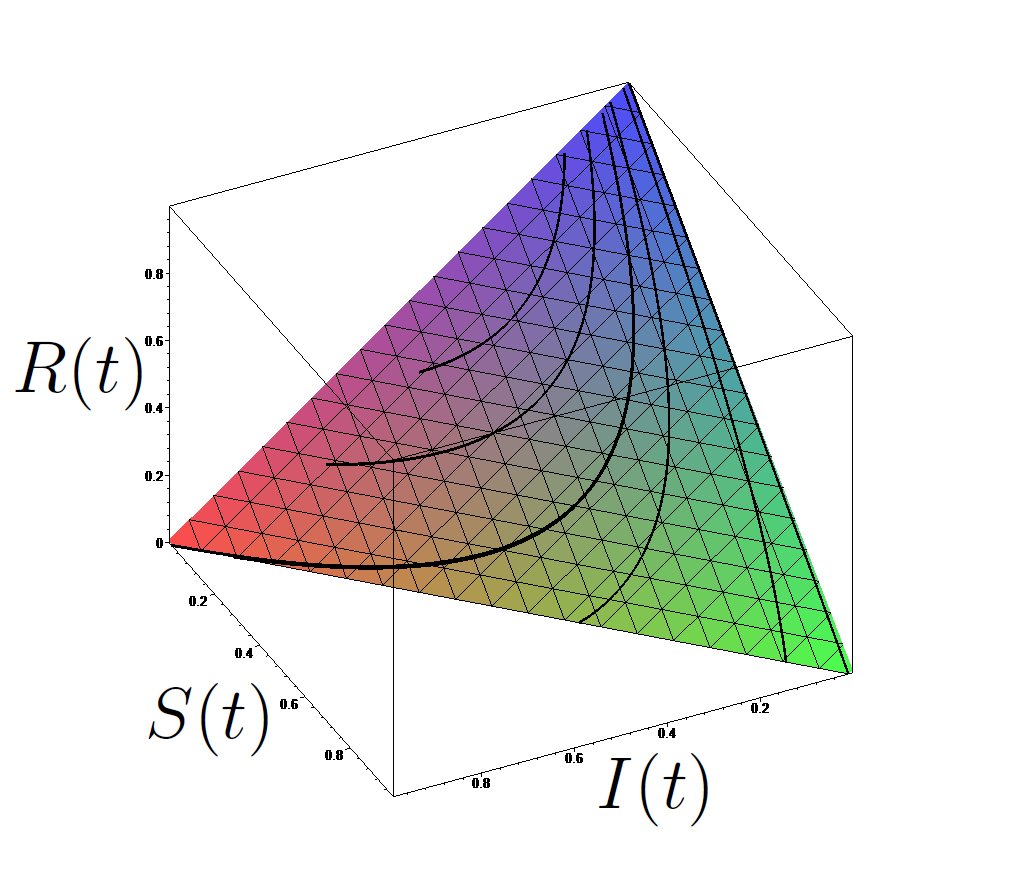}
 \caption{Typical behavior of the trajectories of system \eqref{eq-1} for different initial values.
 Notice that in this model the population of the spreaders tends toward the extinction. The dynamics
 under the constraint $I + S + R =1.$}
 \end{figure}

First of all, we look for the equilibrium points, that is the zeros of the vector field
$$f: {\mathbb R}^3 \to {\mathbb R}^3,\quad f(I,S,R)=(\rho_2  IS,\rho_2  IS-\rho_1  S(S+R), \rho_1  S(S+R)).$$
We also denote by $\nabla{f}$ the corresponding Jacobian matrix
\footnote{We have fixed a minor misprint found in the original paper \cite[Ch. 2, p.3]{Pi-10},
where the variable $S$ is missing in system \eqref{eq-1} and there is
a wrong column in the Jacobian matrix.} which is defined as
$$\nabla{f}(I,S,R)=
\begin{pmatrix}
\frac{\partial I'}{\partial I} & \frac{\partial I'}{\partial S} &  \frac{\partial I'}{\partial R} \\
\frac{\partial S'}{\partial I} & \frac{\partial S'}{\partial S} &  \frac{\partial S'}{\partial R} \\
\frac{\partial R'}{\partial I} & \frac{\partial R'}{\partial S} &  \frac{\partial R'}{\partial R} \\
\end{pmatrix}
\;=\;
\begin{pmatrix}
-\rho_2  & 0 & 0 \\
 \rho_2 S & \rho_2  I + 2\rho_1  S - \rho_1  R & -\rho_1  S\\
0 & 2\rho_1  S + \rho_1  R & \rho_1  S \\
\end{pmatrix}.
$$
We restrict the study of the vector field to the domain
$$D:= \{(I,S,R): I,S,R \geq 0, \, I + S+ R =1\}.$$
By the nature of the constraints defining the domani $D,$ the only possible equilibrium points
are with $S=0$ and then $f(I,0,R)=0$ when $I+R=1$.
As shown by Figures $1$, $2$ and $3$, we have the extinction of the Spreaders when the time
tends to the $+\infty$ since the trajectories tends to equilibrium points
which stay on the hyperplane $I+R=1$.

\begin{figure}[h]
\centering
\includegraphics[width=0.5\textwidth]{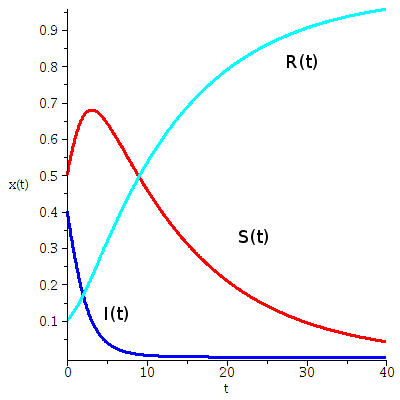}
 \caption{The behavior of the three populations for  $\rho_1=0.1$, $\rho_2=0.9$, $\mu=0.8$
 and initial data $I(0)=0.4$, $S(0)=0.5$, $R(0)=0.1$. The blue, red and cyan lines represent, respectively, the Ignorant,
 Spreaders and Stiflers at the time $t$.
}
\end{figure}

\begin{figure}[h]
\centering
 \includegraphics[width=0.5\textwidth]{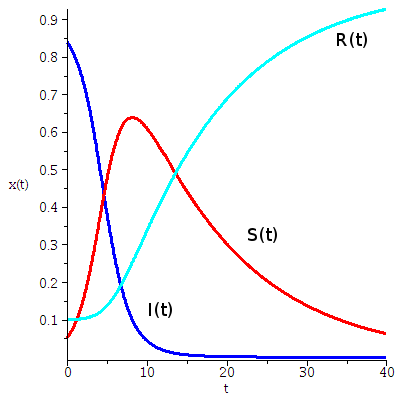}
 \caption{The behavior of the three populations for  $\rho_1=0.1$, $\rho_2=0.9$, $\mu=0.8$
 and initial data $I(0)=0.84$, $S(0)=0.05$, $R(0)=0.1$. The blue, red and cyan lines represent, respectively, the Ignorant,
 Spreaders and Stiflers at the time $t$.
}
 \end{figure}

\newpage

If we compute the Jacobian on these equilibrium points, we obtain:

\begin{equation}
\nabla{f}(I,0,R)=
\begin{pmatrix}
-\rho_2  & 0 & 0 \\
 0 & \rho_2  I - \rho_1  R & 0\\
0 & \rho_1  R & 0 \\
\end{pmatrix}
\end{equation}
The corresponding eigenvalues are:
$\lambda_1 =0, \lambda_2 = -\rho_2 \mu, \lambda_3 = \rho_2 \mu I - \rho_1 \mu R$.
Thus we conclude that $\lambda_2$ is always negative, while
 $\lambda_3 < 0 $ if and only if $(\rho_1+\rho_2)I-\rho_1<0$, that is
$$I<\sigma:=\frac{\rho_1}{\rho_1+\rho_2}.$$
As observed in \cite{Pi-10}, the constant $\sigma$ plays the role of a threshold.
The conclusion in \cite[p.3]{Pi-10} is the following
\begin{itemize}
\item[$(i)$\;]
if $0 < I < \sigma$, the equilibrium point is asymptotically stable;
\item[$(ii)$\;]
if $\sigma < I < 1$, the equilibrium point is unstable.
\end{itemize}
Unfortunately, the conclusion in $(i)$ does not seem completely correct [cf. Figures $4$ and $5$]. 
We recall that an equilibrium point $P$
is \textit{asymptotically stable} if it is stable and there is a neighborhood ${\mathcal U}$ of $P$
such that for each point $z_0$ in ${\mathcal U}$ the solution departing from $z_0$ at the time $t=0$
tends to $P$ as $t\to +\infty.$ We claim that the correct conclusion for $(i)$ would be that for $0 < I < \sigma,$
the equilibrium point is only stable. 

\subsection{First Integral}
A search for stable, non asymptotically, solutions lead us to the search
of a First Integral. By here we show how to proceed with another method that permits to simplify the analysis and
also gives, to our opinion, a better explanation of the results.
Using condition \eqref{cond0} in the normalized form with $N=1$, we can set
$$S(t)=1-I(t)-R(t)$$
and remove the second equation from system \eqref{eq-2}. In this manner, the
original system can be downgraded to a planar system in the two variables $I(t)$ and $R(t)$
that we write as

\begin{equation}\label{eq-3}
\begin{cases}
R'=  \rho_1  (1-I-R)(1-I)\\
I'=-\rho_2  I (1-I-R).\\
\end{cases}
\end{equation}
We denote by  $g(R,I)$ the corresponding vector field related to \eqref{eq-3}.
The analysis of system \eqref{eq-3} will be performed in the set $\Omega$, defined as
$${\Omega}:=\{(R,I): 0 \leq R \leq 1, 0 \leq I \leq 1, R+I \leq 1\}.$$
A simple investigation of the vector field on the boundary of $\Omega$ shows that
on the segment $\{(R,0): 0\leq R \leq 1\},$
we have $R' \geq 0$ and $I' = 0$, while, on the segment $\{(0,I): 0\leq I \leq 1\}.$
we have $R' \geq 0$ and $I'=-\rho_2  I (1-I)$. Furthermore, all the points of the segment
$$\mathscr{S}:=\{(R,I): 0 \leq R \leq 1, 0 \leq I \leq 1, R+I =1\}$$
are equilibrium points. By the uniqueness of the solutions for the initial value problems
associated to \eqref{eq-3}, we conclude that the interior of $\Omega$ is a positively invariant set,
that is $(I(t),R(t))\in \text{int}{\Omega}$ for all $t\geq 0,$ whenever $(I(0),R(0))\in \text{int}{\Omega}$.
The special feature of system \eqref{eq-3} is that there is a \textit{continuum} of equilibrium points for
the equation. Indeed, as previously observed, the set of equilibrium points contained in the domain $\Omega$ is given by
the segment $\mathscr{S}.$ Clearly, such points may be stable or unstable, but they can never be asymptotically stable
(since any neighborhood of an equilibrium point contains infinitely many other equilibria, that is constant solutions).
\\
We consider now the Jacobian matrix $\nabla g$ associated to the two dimensional vector field $g$
and computed at a generic equilibrium point such that $R+I =1.$ A standard computation yields to the following
$$\nabla g(R,I)=\begin{pmatrix}
\rho_1 I-\rho_1  & \rho_1 I-\rho_1\\
\rho_2 I & \rho_2 I
\end{pmatrix}
$$
Clearly, one of the eigenvalues is zero (this is obvious). The other one is give by
$$\tau:=(\rho_1+\rho_2) I-\rho_1.$$
Thus we get the conclusion that if $ 0<I<\sigma$ then equilibrium point  $(R,I)$
is stable, while, if $\sigma<I<1$ the equilibrium point $(R,I)$ is unstable.

In order to provide a more precise description of the global dynamics, we observe that system \eqref{eq-3}
possess a first integral, which can be found via the following steps. First of all, we write equation \eqref{eq-3} as
\begin{equation*}
\begin{cases}
\frac{R'}{\rho_1}= (1-I-R)(1-I)\\
\frac{I'}{I\rho_2}= -(1-I-R),\\
\end{cases}
\end{equation*}
then, we multiply by $1-I)$ the second equation and obtain
\begin{equation*}
\begin{cases}
\frac{R'}{\rho_1}= (1-I-R)(1-I)\\
\frac{I'}{I\rho_2}(1-I)= -(1-I-R)(1-I).\\
\end{cases}
\end{equation*}
Finally, summing up the two equations, we obtain $\frac{I'}{I\rho_2}(1-I)+\frac{R'}{\rho_1}=0$, from which we find that
\begin{equation}\label{eq-4}
\frac{d}{dt} {\mathcal H}(R(t),I(t)) = 0, \quad \text{ for } {\mathcal H}(R,I):=
\frac{R}{\rho_1} + \frac{log(I)}{\rho_2}-\frac{I}{\rho_2}\,.
\end{equation}.
We thus conclude that the solutions of \eqref{eq-3} lie on the level lines of the Hamiltonian function ${\mathcal H},$
that is any solution satisfied the relation
$$\frac{R(t)}{\rho_1} + \frac{log(I(t))}{\rho_2}-\frac{I(t)}{\rho_2} = k = {\mathcal H}(R(0),I(0)), \quad \forall\, t.$$

 \begin{figure}[h]
 \centering
 \includegraphics[width=0.7\textwidth]{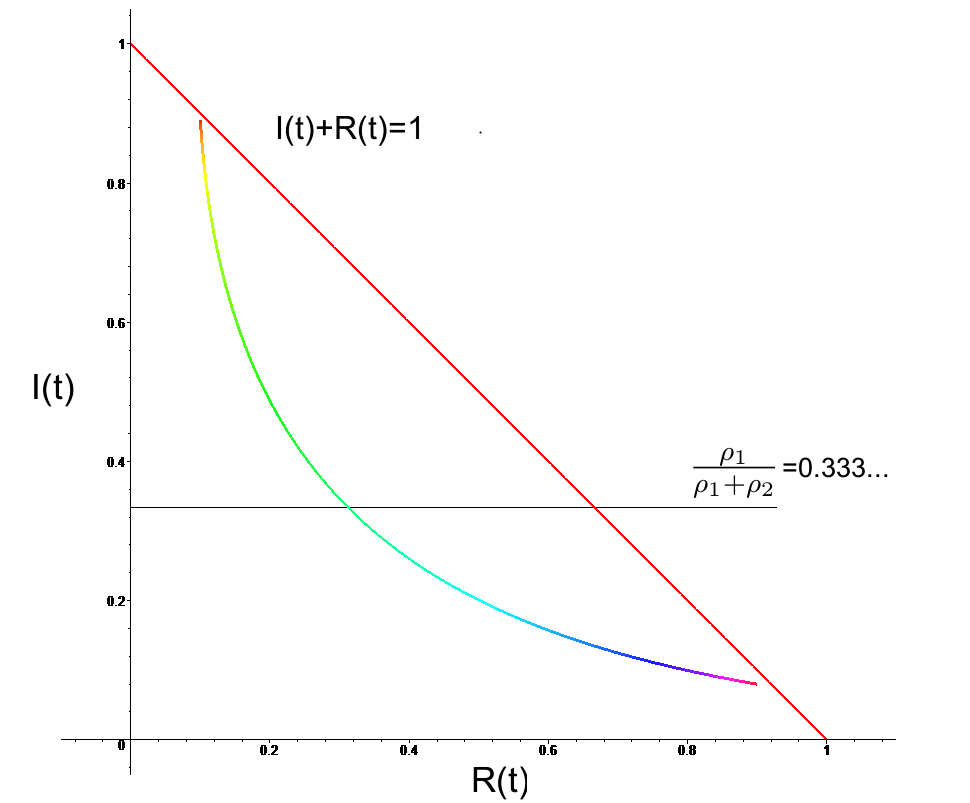}
 \caption{The present figure shows the effect of the threshold value $\sigma$.
 The equilibrium points above the line $I= \sigma$ are of repulsive type. The trajectories
 move from above $\sigma$ and tend asymptotically to some point below the level $\sigma$.
 Lying on the level line of the first integral. The simulation has been performed for $\rho_1=0.4,$ $\rho_2=0.8$
 and $\mu =1.$}
 \end{figure}

 \begin{figure}[h]
 \centering
 \includegraphics[width=0.7\textwidth]{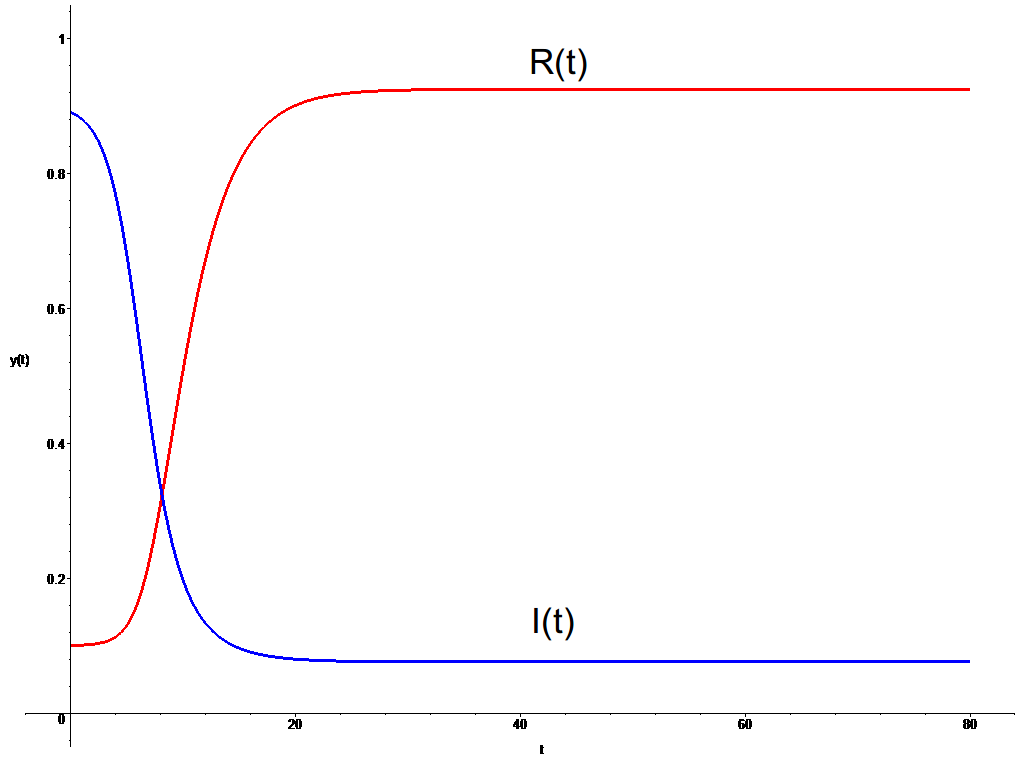}
 \caption{It the same setting of the previous figure, we show the behavior of the solutions for
initial data $I(0)=0.89$ and $R(0)=0.11$. The blue line represents the Ignorant
 at the time $t$ and the red line represents the Stiflers at the time $t$. The simulation shows that
 the two populations quickly stabilize in time.}

\end{figure}

\newpage

\subsection{Other Remarks}
The approach we have described above in the study of Piqueira's model, can be easily adapted to investigate other
rumor transmission models considered by different authors. For instance, we can apply our considerations to a model
by Belen and Pearce in \cite{BePe-04}, where the Authors introduced the differential system:

\begin{equation*}
 \begin{cases}
I'=-IS\\
S'= -S(1-2I) \\
R'=  S(1-I)\\
 \end{cases}
\end{equation*}
With $I(0)=\alpha$, $S(0)=\beta$, $R(0)=\gamma$, $\alpha+\beta+\gamma=1$, $\alpha,\beta >0$ and $\gamma\geq0$.
We observe that is possible remove the equation $R'=  S(1-I)$ and study directly the planar system
\begin{equation*}
\begin{cases}
I'=-IS\\
S'= -S(1-2I). \\
 \end{cases}
\end{equation*}
Passing to the equivalent system
\begin{equation*}
\begin{cases}
\frac{I}{I}'=-S\\
\frac{S'}{S}= (1-2I), \\
 \end{cases}
\end{equation*}
we can easily find a first integral of the form ${\mathcal H}(I,S):= \log(I) - 2I + S.$

We thank prof. Zanolin of the University of Udine for useful discussions about this topic.

\end{document}